\documentclass[12pt]{amsart}

\usepackage{amsthm, amssymb, amsmath}
\usepackage{enumerate}

\usepackage[colorlinks]{hyperref}

\newtheorem{thm}{Theorem}[section]

\newtheorem{lem}[thm]{Lemma}
\newtheorem{prop}[thm]{Proposition}
\newtheorem{cor}[thm]{Corollary}
\newtheorem{rem}[thm]{Remark}
\newtheorem*{nota}{Notation}

\title[Arens regularity of $A_\Phi(G)$ and its dual]{Arens regularity of the Orlicz Fig\`a-Talamanca Herz Algebra}
\author{Arvish Dabra}
\address{Arvish Dabra,\newline\indent Department of Mathematics,\newline\indent Indian Institute of Technology Delhi,\newline\indent New Delhi - 110016, India.}
\email{arvishdabra3@gmail.com}
\author{N. Shravan Kumar}
\address{N. Shravan Kumar,\newline\indent Department of Mathematics,\newline\indent Indian Institute of Technology Delhi,\newline\indent New Delhi - 110016, India.}
\email{shravankumar.nageswaran@gmail.com}

\begin{document}
	
	\begin{abstract}
		Let G be a locally compact group and let $A_\Phi(G)$ be the Orlicz-version of the Fig\`{a}-Talamanca Herz algebra of G associated with a Young function $\Phi.$ We show that if $A_\Phi(G)$ is Arens regular, then $G$ is discrete. We further explore the Arens regularity of $A_\Phi(G)$ when the underlying group $G$ is discrete. In the running, we also show that $A_\Phi(G)$ is finite-dimensional if and only if $G$ is finite. Further, for amenable groups, we show that $A_\Phi(G)$ is reflexive if and only if $G$ is finite, under the assumption that the associated Young function $\Phi$ satisfies the MA-condition.
	\end{abstract}
	
	
	
	\maketitle
	
	\section{Introduction}
	Arens regularity is an important tool to study groups with the help of certain Banach algebras related to it. Richard F. Arens \cite{arens} introduced the first and second Arens products on the double dual $\mathcal{A}^{**}$ of a Banach algebra $\mathcal{A}.$ If both the Arens products coincide, then the Banach algebra $\mathcal{A}$ is said to be Arens regular. Civin and Yood \cite{CY} proved that if $G$ is an infinite locally compact abelian group, then the group algebra $L^1(G)$ cannot be Arens regular. Later, Young \cite{young} proved that for any locally compact group $G,$ $L^1(G)$ is Arens regular if and only if $G$ is finite. In 1964, Eymard \cite{Eym1} introduced and studied the Fourier algebra $A(G),$ for a locally compact group $G$. The above-mentioned result of Young was then generalized to the Fourier algebra $A(G)$ by Lau \cite{Lau1} under the assumption that the group $G$ is amenable.
	
	In 1971, Herz \cite{herz} studied the $p$-version ($1<p<\infty$) of the Fourier algebra. This algebra denoted $A_p(G),$ was later called as Fig\`a-Talamanca Herz algebra. For $p=2,$ $A_p(G)$ coincides with the Fourier algebra $A(G).$ In \cite{For1}, Forrest showed that the Arens regularity of the algebra $A_p(G)$ implies that the group $G$ is discrete. Later, he proved that for an abelian group $G, A_p(G)$ is Arens regular if and only if $G$ is finite \cite{For2}. Lau and \"{U}lger \cite{LU} proved an equivalent condition for the Arens regularity of $A_p(G)$ in terms of the spaces $PM_p(G)$ and $PF_p(G).$
	
	It is well known that Orlicz spaces are the natural generalization of the classical $L^p$-spaces. The Orlicz-version of the Fig\`a-Talamanca Herz algebra, denoted $A_\Phi(G),$ was recently introduced and studied by the authors of \cite{RLSK1}. Aghababa and Akbarbaglu \cite{AA} independently introduced and studied the same algebra in 2020. In this paper, we shall study the structure of the second dual $A_\Phi(G)^{\ast\ast}$ as a Banach algebra concerning the two Arens products, introduced by Arens \cite{arens}. In particular, we shall show that if $A_\Phi(G)$ is Arens regular, then G is discrete. This is Theorem \ref{thm2} of this paper.
	
	It is a classical problem of Functional analysis to characterize Banach spaces that are reflexive. Under the assumption that the group $G$ is amenable, Granirer \cite{Gra1} showed that the $A_p(G)$ algebras are reflexive if and only if $G$ is finite. Corollary \ref{Ref_Chac} of the present paper gives a similar result for the $A_\Phi(G)$ algebras. We also characterize the property of being finite-dimensional for the algebras $A_\Phi(G).$ In fact, we shall show that $A_\Phi(G)$ is finite-dimensional if and only if $G$ is finite. This is Theorem \ref{thm8} of the present paper.
	
	In section 4, we provide an equivalent condition in terms of certain spaces associated with $A_\Phi(G).$ The approach in the proof is inspired by Lau and \"{U}lger's methods in \cite{LU}. In section 5, we again study the Arens regularity of $A_\Phi(G),$ assuming that the group $G$ is discrete. The authors of \cite{RLSK4} studied the space of pseudofunctions, denoted $PF_\Psi(G)$ and its dual $W_\Phi(G).$ Section 6 concludes this paper with a study of the space $PF_\Psi(G).$ We shall prove that the algebra $PF_\Psi(G)$ is Arens regular for a compact group $G$ and that its dual $W_\Phi(G)$ exhibits the Radon-Nikod\'ym property. 
	
	To prepare for the sequel, we shall start with some preliminaries.
	
	\section{Preliminaries}
	
	A symmetric convex function $\Phi: \mathbb{R} \to [0,\infty]$ which satisfies $\Phi(0) = 0$ and $\lim\limits_{x \to \infty} \Phi(x) = +\infty$ is referred to as Young function. With each Young function $\Phi,$ another convex function $\Psi: \mathbb{R} \to [0,\infty]$ can be associated, given by $$\Psi(y) = \sup\{x \, |y| - \Phi(x): x \geq 0\}, \hspace{1cm} y \in \mathbb{R}.$$ Then the function $\Psi$ is also a Young function and is termed as the complementary function to $\Phi.$ Consequently, the pair $(\Phi,\Psi)$ is referred to as complementary pair of Young functions. It is easy to verify that if $(\Phi,\Psi)$ is a complementary pair of Young functions, then so is the pair $(\Psi,\Phi).$
	
	A continuous Young function $\Phi$ is termed as N-function if it satisfies $\Phi(x) = 0$ iff $x = 0$ and $\lim\limits_{x \to 0} \frac{\Phi(x)}{x} = 0, \lim\limits_{x \to \infty} \frac{\Phi(x)}{x} = +\infty$ while $\Phi(\mathbb{R}) \subseteq \mathbb{R}^+.$ For each $p \in (1,\infty),$ consider the Young function $\Phi_p(x) = \frac{|x|^p}{p}.$ It is easy to check that $\Phi_p$ is an N-function with $\Psi_q(y) = \frac{|y|^q}{q}, \, (p^{-1}+q^{-1}=1),$ as its complementary function. It should be noted that $\Phi_1$ is only a Young function, not an N-function. Consider the functions $\Phi_0(x) = e^{|x|}-|x|-1$ and $\Psi_0(y) = (1+|y|)\log(1+|y|)-|y|.$ Then, $(\Phi_0,\Psi_0)$ is another example of a complementary pair of Young functions.
	
	Let $G$ be a locally compact group with a left Haar measure $dx.$ A Young function $\Phi$ is said to satisfy the $\Delta_2$-condition, denoted $\Phi \in \Delta_2,$ if there exists a constant $M > 0$ and $x_0 > 0$ such that $\Phi(2x) \leq M \, \Phi(x)$ whenever $x \geq x_0$ if $G$ is compact and the same inequality holds with $x_0 = 0$ if $G$ is non compact. The pair $(\Phi,\Psi)$ of complementary Young functions is said to satisfy the $\Delta_2$-condition if both $\Phi$ and $\Psi$ satisfies the $\Delta_2$-condition. As an example, the Young function $\Psi_0(y) = (1+|y|)\log(1+|y|)-|y| \in \Delta_2,$ however, its complementary function $\Phi_0(x) = e^{|x|}-|x|-1 \notin \Delta_2.$ It is easy to verify that the family of Young functions $\Phi_{\alpha,p}: y \mapsto \alpha \, |y|^p$ for $p \geq 1$ and $\alpha > 0,$ satisfies the $\Delta_2$-condition.
	
	The Orlicz space, denoted $L^\Phi(G),$ is defined as $$L^\Phi(G) = \bigg\{f : G \to \mathbb{C}: f \, \text{is measurable and} \, \int\limits_G \Phi(\beta \, |f|)dx < \infty \, \text{for some} \, \beta > 0\bigg\}.$$ The Orlicz space $L^\Phi(G)$ is a Banach space when equipped with the norm $$N_\Phi(f) = \inf\bigg\{k > 0: \int\limits_G \Phi\bigg(\frac{|f|}{k}\bigg)dx \leq 1\bigg\}.$$ The aforementioned norm is called the Gauge norm or Luxemburg norm. If $(\Phi,\Psi)$ is a complementary pair of Young functions, then there is another norm on $L^\Phi(G),$ given by $$\|f\|_\Phi = \sup\bigg\{\int\limits_G |fg|: g \in L^\Psi(G) \, \text{and} \int\limits_G \Psi(|g|)dx \leq 1\bigg\}.$$ This norm is called the Orlicz norm. In fact, these two norms are equivalent, and for any $f \in L^\Phi(G),$ $$N_\Phi(f) \leq \|f\|_\Phi \leq 2 N_\Phi(f).$$
	
	Let $\mathbf{C}_c(G)$ denote the space of all continuous functions on $G$ with compact support. If a Young function $\Phi$ satisfies the $\Delta_2$-condition, then $\mathbf{C}_c(G)$ is norm dense in $L^\Phi(G).$ Furthermore, if the complementary function $\Psi$ is continuous and satisfies $\Psi(x) = 0$ iff $x = 0,$ then the dual of $(L^\Phi(G),N_\Phi(\cdot))$ is isometrically isomorphic to $(L^\Psi(G),\|\cdot\|_\Psi)$ \cite[Corollary 9, Pg. 111]{raoren}. In particular, if the complementary Young pair $(\Phi,\Psi)$ satisfies the $\Delta_2$-condition, then the Banach space $L^\Phi(G)$ is reflexive \cite[Theorem 10, Pg. 112]{raoren}.
	
	Let $G$ be a compact group. Then, for any $f \in L^1(G),$ the map $g \mapsto f \ast g$ (also, $g \mapsto g \ast f)$ from $L^\Phi(G)$ to $L^\Phi(G)$ is compact. Further, $L^\Phi(G)$ has an approximate identity that is bounded in $\|\cdot\|_1$.
	
	A Young function $\Phi$ is said to satisfy the \textit{Milnes-Akimoni$\check{c}$} condition (in short, MA-condition) if for each $\epsilon > 0,$ there exists $\alpha_\epsilon > 1$ and an $x_1(\epsilon) \geq 0$ such that $$\Phi'((1+\epsilon)x) \geq \alpha_\epsilon \, \Phi'(x), \hspace{1cm} x \geq x_1(\epsilon).$$ The Young function $\Phi(x) = e^x - x - 1, \, 0 \leq x \in \mathbb{R},$ is an example of an N-function that satisfies the MA-condition. Although its complementary function $\Psi(y) = (1+y)\log(1+y)-y, \, 0 \leq y \in \mathbb{R},$ is also an N-function, but it does not satisfy the MA-condition. The family of Young functions $\Phi_\alpha(x) = |x|^\alpha (1 + |\log|x||),$ $\alpha > 1,$ satisfies the MA-condition. Moreover, the family $(\Phi_\alpha,\Psi_\alpha)$ of complementary pair of Young functions also satisfies the $\Delta_2$-condition.
	
	For more details on Orlicz spaces, the reader can refer to \cite{raoren}.
	
	Let $(\Phi,\Psi)$ be a complementary Young pair satisfying the $\Delta_2$-condition. If $h: G \to \mathbb{C}$ is any function then $\check{h}$ is defined as $\check{h}(x) = h(x^{-1})$ for $x \in G$. Let $$A_\Phi(G)=\left\{u =\sum\limits_{n \in \mathbb{N}} f_n \ast \check{g_n}:
	\begin{array}{c}
		\{f_n\} \subseteq L^\Phi(G), \{g_n\} \subseteq L^\Psi(G),  \\
		\sum\limits_{n \in \mathbb{N}} N_\Phi(f_n) \|g_n\|_\Psi<\infty
	\end{array} \right\}.$$ 
	Note that if $u \in A_\Phi(G),$ then $u \in \mathbf{C}_0(G).$ For $u \in A_\Phi(G),$ define $\|u\|_{A_\Phi}$ as $$\|u\|_{A_\Phi} = \inf\bigg\{\sum\limits_{n \in \mathbb{N}} N_\Phi(f_n) \|g_n\|_\Psi: u = \sum\limits_{n \in \mathbb{N}} f_n \ast \check{g_n}\bigg\}.$$ Equipped with this norm and pointwise addition and multiplication, the space $A_\Phi(G)$  becomes a commutative Banach algebra. This Banach algebra $A_\Phi(G)$ is called the Orlicz Fig\`a-Talamanca Herz algebra. In fact, $A_\Phi(G)$ is semi-simple and regular, with the spectrum being homeomorphic to $G.$
	
	Let $B_\Phi(G)$ denote the space $\{u \in \mathbf{C}(G): uv \in A_\Phi(G) \, \forall \, v \in A_\Phi(G)\}$. Each $u \in B_\Phi(G)$ defines a bounded linear map $M_u: A_\Phi(G) \to A_\Phi(G),$ given by $M_u(v) = uv.$ Then $B_\Phi(G)$ becomes a commutative Banach algebra when equipped with the operator norm and pointwise addition and multiplication. It is easy to see that $A_\Phi(G) \subseteq B_\Phi(G).$
	
	Let $\mathbf{B}(L^\Phi(G))$ be the space of all bounded linear operators on $L^\Phi(G)$ equipped with the operator norm. For a bounded complex Radon measure $\mu$ on $G$ and $f \in L^\Phi(G),$ define $T_\mu: L^\Phi(G) \to L^\Phi(G)$ by $T_\mu(f) = \mu \ast f.$ It is clear that $T_\mu \in \mathbf{B}(L^\Phi(G)).$ Let $PM_\Phi(G)$ denote the closure of $\{T_\mu: \mu \in M(G)\}$ in $\mathbf{B}(L^\Phi(G))$ with respect to the ultra-weak topology (weak$^*$-topology). It is proved in \cite{RLSK1} that for a locally compact group $G,$ the dual of $A_\Phi(G)$ is isometrically isomorphic to $PM_\Psi(G).$ Let $CV_\Psi(G)$ denote the space of all bounded linear operators $T: L^\Psi(G) \to L^\Psi(G)$ such that $T(f\ast g) = T(f) \ast g$ for all $f \in L^\Psi(G)$ and $g \in L^1(G).$ It is clear that $PM_\Psi(G)$ is contained in $CV_\Psi(G)$ and in fact, this inclusion map is an $A_\Phi(G)$-module map.  Let $PF_\Phi(G)$ denote the norm closure of $\{T_f: f \in L^1(G)\}$ inside $\mathbf{B}(L^\Phi(G)).$ The dual of $PF_\Phi(G)$ is denoted by $W_\Psi(G).$
	
	A topological invariant mean on $PM_\Psi(G)$ is a functional $m \in PM_\Psi(G)^*$ satisfying $\|m\| = 1 = m(I)$ and $u \cdot m = u(e) m \,  \forall \, u \in A_\Phi(G),$ i.e., $$\langle T,u\cdot m\rangle = \langle u \cdot T,m\rangle = u(e) \langle T,m \rangle \, \, \forall \, u \in A_\Phi(G) \, \, \text{and} \, \, \forall \, \, T \in PM_\Psi(G).$$ The set of all topological invariant means on $PM_\Psi(G)$ is denoted by $TIM_\Psi(\widehat{G}).$
	
	For further details on $A_\Phi(G), PM_\Psi(G), PF_\Psi(G)$ and $TIM_\Psi(\widehat{G}),$ one can refer to the papers \cite{AA}, \cite{RLSK1} and \cite{RLSK4}.
	
	Let $\mathcal{A}$ be a Banach algebra. Then, the double dual $\mathcal{A}^{**}$ can be given two multiplications which extend the multiplication of $\mathcal{A}$ and for which $\mathcal{A}^{**}$ becomes a Banach algebra. Arens introduced these two products for normed algebras in \cite{arens}. These two products are defined as follows:
	\begin{enumerate}[1)]
		\item $\odot$ denotes the first Arens product.
		\begin{enumerate}[(i)]
			\item $\langle u \cdot T,v\rangle = \langle T,vu\rangle \, \text{for every} \, u,v \in \mathcal{A}, T \in \mathcal{A}^*.$
			\item $\langle T \odot m,u\rangle = \langle m,u \cdot T\rangle \, \text{for every} \, u \in \mathcal{A}, T \in \mathcal{A}^*, m \in \mathcal{A}^{**}.$
			\item $\langle m_1 \odot m_2,T\rangle = \langle m_2,T \odot m_1\rangle \, \text{for every} \, T \in \mathcal{A}^*, m_1, m_2 \in \mathcal{A}^{**}.$
		\end{enumerate}
		\item $\square$ denotes the second Arens product.
		\begin{enumerate}[(i)]
			\item $\langle T \, \square \, u,v\rangle = \langle T,uv\rangle \, \text{for every} \, u,v \in \mathcal{A}, T \in \mathcal{A}^*.$
			\item $\langle m \, \square \, T,u\rangle = \langle m,T \, \square \, u\rangle \, \text{for every} \, u \in \mathcal{A}, T \in \mathcal{A}^*, m \in \mathcal{A}^{**}.$
			\item $\langle m_1 \, \square \, m_2,T\rangle = \langle m_1,m_2 \, \square \, T\rangle \, \text{for every} \, T \in \mathcal{A}^*, m_1, m_2 \in \mathcal{A}^{**}.$
		\end{enumerate}
	\end{enumerate}
	In general, $m_1 \odot m_2$ is not equal to $m_1 \, \square \, m_2$ for all $m_1, m_2 \in \mathcal{A}^{**}.$ If this holds for all $m_1, m_2 \in \mathcal{A}^{**},$ then $\mathcal{A}$ is said to be Arens regular.
	
	Let $\mathcal{A}$ be a commutative Banach algebra. Then, $u \cdot T = T \, \square \, u.$ Moreover, $\mathcal{A}$ is Arens regular iff $\mathcal{A}^{**}$ becomes a commutative Banach algebra concerning either, and hence both, of the Arens products \cite[Proposition 1]{DH}.
	
	A functional $T \in \mathcal{A}^*$ is said to be an almost periodic functional on $\mathcal{A}$ if the mapping $a \mapsto a \cdot T$ from $\mathcal{A}$ to $\mathcal{A}^*$ is a compact operator. Similarly, a weakly almost periodic functional on $\mathcal{A}$ is a functional $T \in \mathcal{A}^*$ such that the mapping $a \mapsto a \cdot T$ from $\mathcal{A}$ to $\mathcal{A}^*$ is a weakly compact operator. For $\mathcal{A} = A_\Phi(G),$ the space of all almost periodic and weakly almost periodic functionals are denoted $AP_\Psi(\widehat{G})$ and $WAP_\Psi(\widehat{G}),$ respectively. A functional $T \in \mathcal{A}^*$ is called uniformly continuous if $T$ is in the norm closure of $\text{span}\{u \cdot \widehat{T}: u \in \mathcal{A} \, \text{and} \, \widehat{T} \in \mathcal{A}^*\}.$ When $\mathcal{A} = A_\Phi(G),$ the space of uniformly continuous functionals is denoted by $UCB_\Psi(\widehat{G}).$
	
	We refer the readers to \cite{RLSK3} for more details on $AP_\Psi(\widehat{G}),$ $WAP_\Psi(\widehat{G})$ and $UCB_\Psi(\widehat{G}).$ 
	
	A Banach space $X$ has the Radon-Nikod\'ym property (RNP) if each closed convex subset $D$ of $X$ is dentable, i.e., for any $\epsilon > 0,$ there exists an $x \in D$ such that $x \notin \overline{co}(D \setminus B_\epsilon(x)),$ where $B_\epsilon(x) = \{y \in X: \|x-y\| < \epsilon\}.$ \cite{Bou}
	
	Throughout this paper, $G$ shall denote a locally compact group with a fixed left Haar measure, denoted $dx$ and sometimes as $m.$ Further, $\Phi$ denotes a Young function with $\Psi$ as its complementary function, and the pair $(\Phi,\Psi)$ satisfies the $\Delta_2$-condition.
	
	\section{Reflexivity and finite-dimensionality of $A_\Phi(G)$}
	In this section, we characterize when $A_\Phi(G)$ is finite-dimensional and reflexive, respectively. We shall see that both properties individually characterize finite groups. To achieve this, we first characterize reflexive $A_\Phi(G)$-submodules of $PM_\Psi(G).$ The results of this section are motivated from \cite{Gra1, Gra2}.
	
	We shall begin with the following lemma.
	
	\begin{lem}\label{lemma}
		Let $h \in PF_\Psi(G)^*.$ For $\mu \in M(G),$ define $$\langle h,\mu \rangle = \int\limits_G h(x) d\mu(x).$$ Then $h$ defines a bounded linear functional on $(M(G),\| \cdot \|_{PM_\Psi}).$
	\end{lem}
	\begin{rem}
		Note that $M(G)$ sits inside $PM_\Psi(G)$ via the following identification: $\mu \mapsto M_\mu,$ where $M_\mu(v) = \int\limits_G v \, d\mu$ for $v \in A_\Phi(G)$ or $M_\mu(g) = \mu \ast g$ for $g \in L^\Psi(G).$
	\end{rem}
	\begin{rem}\label{rem:neat}
		Note that for any $f \in L^1(G),$ the measure $\mu_f$ given by $$\mu_f(E) = \int\limits_G \chi_E(x) f(x) \, dx$$ belongs to $M(G).$ Thus, for any $T \in PF_\Psi(G)^* = W_\Phi(G),$ we can define $\varphi_T: L^1(G) \to \mathbb{C}$ as $\varphi_T(f) = T(M_{\mu_f}).$ Then
		\begin{align*}
			|\varphi_T(f)| &= |T(M_{\mu_f})| \leq \|T\|_{PF_\Psi(G)^*} \, \|M_{\mu_f}\|_{PF_\Psi(G)} \leq \|T\|_{W_\Phi(G)} \, \|f\|_1,
		\end{align*}
		i.e., $\varphi_T \in L^1(G)^*.$ Therefore, there exists $h_T \in L^\infty(G)$ such that $$\varphi_T(f) = \int\limits_G f(x) h_T(x) dx.$$ Further, $$\bigg| \int\limits_G f h_T \bigg| \leq \|h_T\|_{W_\Phi(G)} \, \|f\|_{PF_\Psi(G)}.$$ This is how we identify any element of $W_\Phi(G)$ as a function.
	\end{rem}
	\begin{proof}[Proof of Lemma \ref{lemma}]
		Let $\{V_\alpha\}_{\alpha \in \Lambda}$ be a neighbourhood base of $e$ such that for each $\alpha \in \Lambda,$ $V_\alpha$ has compact closure. Now, for each $\alpha \in \Lambda,$ choose a continuous function $k_\alpha$ such that $k_\alpha \geq 0,$ \text{supp}$(k_\alpha) \subseteq V_\alpha$ and $\int\limits_G k_\alpha = 1.$ Fix $\mu \in M(G).$ Define $f_\alpha: G \to \mathbb{C}$ as $f_\alpha = \mu \ast k_\alpha.$ It is clear that $f_\alpha \in L^1(G)$ and the net $\int\limits_G f_\alpha \, g$ converges to $\int\limits_G g \, d\mu$ for each continuous function $g$ on $G$ \cite[Pg. 190]{Eym1}. By Remark \ref{rem:neat} and \cite[Corollary 3.2]{RLSK4}, we have $$W_\Phi(G) \subseteq B_\Phi(G) \subseteq \mathbf{C}(G).$$ Thus, for $h \in PF_\Psi(G)^* = W_\Phi(G),$
		\begin{align*}
			\bigg | \int\limits_G h(x) d\mu(x) \bigg | &= \underset{\alpha}{\lim} \, \bigg | \int\limits_G f_\alpha(x) h(x) dx \bigg | \leq \underset{\alpha}{\lim} \, \|h\|_{W_\Phi(G)} \, \|f_\alpha\|_{PF_\Psi(G)}\\
			&\leq \|h\|_{W_\Phi(G)} \, \|\mu\|_{PM_\Psi(G)} \, \underset{\alpha}{\lim} \, \|k_\alpha\|_1 = \|h\|_{W_\Phi(G)} \, \|\mu\|_{PM_\Psi(G)},
		\end{align*}
		showing that $h$ defines a bounded linear functional on $(M(G),\| \cdot \|_{PM_\Psi}).$
	\end{proof}
	\begin{nota}
		For $x \in G,$ we write $\lambda_\Psi(x)$ for the operator $T_{\delta_x}$ on $L^\Psi(G).$
	\end{nota}
	\begin{prop}\label{propo}
		Let $G$ be an amenable group and let $\Phi$ satisfy the MA-condition. Let $S \subseteq G$ be such that the set $\{\lambda_\Psi(x): x\in S\}$ is weakly relatively compact in $PM_\Psi(G).$ Then S is finite.
	\end{prop}
	\begin{proof}
		Let $C$ denote the weak closure of $\{\lambda_\Psi(x): x\in S\}$ inside $PM_\Psi(G).$ By assumption, $C$ is weakly compact. Suppose to the contrary that $S$ is infinite. Then, by Eberlein-Smulian theorem, there exists a sequence $\{x_n\} \subseteq S$ such that each element of the sequence is distinct and $\lambda_\Psi(x_n) \to T$ weakly for some $T \in PM_\Psi(G).$ Using the fact that the weak{$^*$}-topology is weaker than the weak topology, we obtain that $\lambda_\Psi(x_n) \to T$ in the weak{$^*$}-topology. As each $\lambda_\Psi(x_n)$ serves as a multiplicative linear functional on $A_\Phi(G),$ either $T$ is also a multiplicative linear functional on $A_\Phi(G)$ or $T = 0.$ Thus, by \cite[Corollary 3.8]{RLSK1}, either $T = \lambda_\Psi(x_0)$ for some $x_0 \in G$ or $T = 0.$
		
		Now, let $m$ be a topological invariant mean on $PM_\Psi(G)$ (such a $m$ exists by \cite[Corollary 6.2]{RLSK1}). By \cite[Theorem 3.4]{RLSK3}, $m$ can be chosen such that, for $\mu \in M(G), \, m(\mu) = \mu(\{e\}),$ where $e$ denotes the identity of the group $G.$ For an arbitrary $x \in G,$ let $m_x$ be the invariant mean on $PM_\Psi(G)$ given by $\langle T,m_x \rangle = \langle \lambda_\Psi(x) T,m \rangle.$ Note that $\lambda_\Psi(x) \, \lambda_\Psi(y) = \lambda_\Psi(xy)$ for all $x,y \in G$ and hence $$\langle \lambda_\Psi(y),m_x \rangle = \langle \lambda_\Psi(xy),m \rangle = 
		\begin{cases} 
			1 & \text{if} \, \, xy = e \\
			0 & \text{otherwise}. 
		\end{cases}$$
		Thus, if $T = \lambda_\Psi(x_0),$ then $$\langle \lambda_\Psi(x_n),m_{{x_0}^{-1}} \rangle = 
		\begin{cases} 
			1 & \text{if} \, \, x_n = x_0 \, (\text{at most one such n is possible}) \\
			0 & \text{otherwise} 
		\end{cases}$$
		and therefore, $\lambda_\Psi(x_n) \nrightarrow \lambda_\Psi(x_0)$ weakly. On the other hand, if $T = 0,$ using the fact that $G$ is amenable and by \cite[Corollary 3.3]{RLSK4}, $1 \in B_\Phi(G) \cong W_\Phi(G) = PF_\Psi(G)^*.$ By Lemma \ref{lemma}, the mapping $\mu \mapsto \langle 1,\mu \rangle$ is a bounded linear functional on $(M(G),\| \cdot \|_{PM_\Psi}).$ Therefore, $1 = \langle 1,\lambda_\Psi(x_n) \rangle \to \langle 1,T \rangle = 0,$ which is absurd.
		
		Thus, we conclude that $S$ is finite.
	\end{proof}
	\begin{nota}
		For any set $X, \, \#X$ shall denote the cardinality of the set $X.$
	\end{nota}
	In the next theorem, we characterize reflexive $A_\Phi(G)$-submodules of $PM_\Psi(G).$
	\begin{thm}\label{thm5}
		Suppose that $G$ is amenable, $\Phi$ satisfies the MA-condition and $N$ is a norm closed $A_\Phi(G)$-submodule of $PM_\Psi(G).$ If $N$ is reflexive, then there exists a finite set $F \subseteq G$ such that $N = span\{\lambda_\Psi(x): x \in F\}.$
	\end{thm}
	\begin{proof}
		Let $B_{N,1}$ denote the closed unit ball of $N.$ By assumption, $N$ is reflexive, and therefore $B_{N,1}$ is weakly compact. Hence, $B_{N,1}$ is $\sigma(PM_\Psi(G), PM_\Psi(G)^\ast)$-compact and this implies that it is w*-compact in $PM_\Psi(G).$ Now, by imitating the proof of \cite[Theorem 1.3]{Gra2}, one can show that $N$ is w*-closed in $PM_\Psi(G).$ Let $F=\{x\in G:\lambda_\Psi(x)\in N\}.$ Then, the set $F_N=\{\lambda_\Psi(x):x\in F\}\subseteq B_{N,1}.$ Since $B_{N,1}$ is weakly compact, $F_N$ is relatively weakly compact and therefore, by Proposition \ref{propo}, $F_N,$ and hence $F,$ is finite. 
		
		We now claim that $N = \text{span}\{\lambda_\Psi(x): x \in F \}.$ In fact, it is enough to show that $N \subseteq \text{span}\{\lambda_\Psi(x): x \in F \}.$ So, let $T \in N.$ Then, for any $x \in \text{supp}(T),$ $\lambda_\Psi(x)$ is the w$^*$-limit of operators of the form $u.T$ with $u \in A_\Phi(G).$ As $N$ is a w$^*$-closed $A_\Phi(G)$-submodule of $PM_\Psi(G),$ we have that $\lambda_\Psi(x) \in N,$ which implies that $x \in F.$ Hence, $\text{supp}(T) \subseteq F.$ Let $\text{supp}(T) = \{ x_k: 1 \leq k \leq l \},$ where $l \leq \#(F).$
		
		Let $V$ be a compact symmetric neighbourhood of $e$ such that $x_k VVV \, \cap \, x_j VVV = \emptyset$ for $k \neq j.$ Let $u = \frac{1}{|V|} (\chi_V \ast \chi_{VV}).$ Then $u \in A_\Phi(G), \, u(x) = 1$ on $V$ and $u = 0$ on $(VVV)^c.$ For $1 \leq i \leq l,$ let $u_i = u_{{x_i}^{-1}}$ where $u_{y}(x) = u(yx).$ Since $\text{supp}(u_i) \subseteq x_i VVV$, we have $\text{supp}(u_i \, . \, T) \subseteq \text{supp}(u_i) \cap \text{supp}(T) = \{x_i\}.$ By \cite[Theorem 3.6]{RLSK1}, singletons are sets of spectral synthesis for $A_\Phi(G)$ and thus $u_i \, . \, T = \alpha_i \lambda_\Psi(x_i)$ for some scalars $\alpha_i.$ Thus $$\bigg (\sum_{1}^{l} u_i \bigg ). \, T  = \sum_{1}^{l} \alpha_i \lambda_\Psi(x_i).$$ Further, if we take $v = \sum\limits_{1}^{l} u_i,$ then $v \in A_\Phi(G)$ and $v \equiv 1$ on some open neighbourhood containing $\text{supp}(T).$ Also, $$T = v \, . \, T = \sum_{1}^{l} \alpha_i \lambda_\Psi(x_i),$$ showing that $T \in \text{span}\{\lambda_\Psi(x): x \in F\}.$
	\end{proof}
	We now proceed towards the main results of this section. We shall begin by characterizing the reflexivity of the quotient space $A_\Phi(G)/I,$ where $I$ is a closed ideal.
	\begin{cor}\label{cor6}
		Let $I$ be a closed ideal in $A_\Phi(G).$ Suppose that $G$ is amenable and that $\Phi$ satisfies the MA-condition. Then $A_\Phi(G)/I$ is reflexive if and only if it is finite-dimensional.
	\end{cor}
	\begin{proof}
		Let $N = (A_\Phi(G)/I)^*.$ Then, by a standard result of functional analysis, $N = \{T \in PM_\Psi(G) : T|_I = 0 \} = I^\perp.$ It is now clear that $N$ is w$^*$-closed. Since $I$ is an ideal, it follows that $N$ is an $A_\Phi(G)$-submodule of $PM_\Psi(G).$ Thus, $A_\Phi(G)/I$ is reflexive if and only if $N$ is reflexive, which by Theorem \ref{thm5}, is equivalent to saying that $N$ is finite-dimensional or equivalently, $A_\Phi(G)/I$ is finite-dimensional.
	\end{proof}
	The next result gives the existence of a function in $A_\Phi(G)$ with certain properties.
	\begin{lem}\label{lem7}
		Let $K \subseteq U \subseteq G$ where $U$ is open and $K$ is compact. Then there exists $u \in A_\Phi(G) \cap \mathbf{C}_c(G)$ such that
		\begin{enumerate}[(i)]
			\item $0 \leq u(x) \leq 1 \, \forall \, x \in G,$
			\item $u \equiv 1$ on $K,$
			\item $\text{supp}(u) \subseteq U.$
		\end{enumerate}
	\end{lem}
	\begin{proof}
		Choose a compact neighbourhood $V$ of $e$ in $G$ such that $VV^{-1}K \subseteq U$ and let $u = \frac{1}{|V|}(\chi_V \ast \chi_{K^{-1}V}) \in A_\Phi(G).$ Then, the function $u$ fulfills the requirements of the lemma.
	\end{proof}
	The following theorem provides the necessary and sufficient condition for the Banach algebra $A_\Phi(G)$ to be finite-dimensional.
	\begin{thm}\label{thm8}
		The Banach algebra $A_\Phi(G)$ is finite-dimensional if and only if $G$ is finite.
	\end{thm}
	\begin{proof}
		If $G$ is finite, it is clear that $A_\Phi(G)$ is finite-dimensional. We shall now prove the converse. The proof is by contradiction. Suppose that $A_\Phi(G)$ is finite-dimensional but $G$ is infinite. Let $\{x_n\} \subseteq G$ be an infinite sequence consisting of distinct elements. Since $x_1 \neq x_2,$ there exists disjoint open sets $U_1$ and $U_2$ containing $x_1$ and $x_2,$ respectively. By Lemma \ref{lem7}, there exists $u_1 \in A_\Phi(G)$ such that $u_1(x_1) = 1$ and $\text{supp}(u_1) \subseteq U_1.$ As $U_1 \cap U_2 = \emptyset, u_1(x_2) = 0.$ Now, $\{x_1,x_2\}$ is compact and therefore, by using the fact that a locally compact Hausdorff space is regular, there exists disjoint open sets $U_{12}$ and $U_3$ such that $\{x_1,x_2\} \subseteq U_{12}$ and $x_3 \in U_3.$ Once again by Lemma \ref{lem7}, there exists $u_2 \in A_\Phi(G)$ such that $u_2(x_1) = u_2(x_2) = 1$ and $u_3(x_3) = 0.$ Proceeding in this way, we can obtain a sequence $\{u_n\} \subseteq A_\Phi(G)$ such that $u_n(x_i) = 1$ for all $1 \leq i \leq n$ and $u_n(x_{n+1}) = 0.$ In particular, we have a linearly independent infinite subset of $A_\Phi(G),$ contradicting the finite-dimensionality of $A_\Phi(G).$ Therefore, $G$ must be finite.
	\end{proof}
	As a consequence of the above theorem, we have the following corollary, which characterizes the reflexivity of the algebra $A_\Phi(H)$ when $H$ is an open subgroup of $G.$
	\begin{cor}\label{cor9}
		Let $H$ be an open subgroup of an amenable group $G.$ Suppose that $\Phi$ satisfies the MA-condition. Then $A_\Phi(H)$ is reflexive if and only if $H$ is finite.
	\end{cor}
	\begin{proof}
		Let $H \subseteq G$ be an open subgroup, and let $I(H) = \{u \in A_\Phi(G): u \equiv 0 \, \text{on} \, H\}.$ By \cite[Lemma 4.1]{RLSK1}, $A_\Phi(G)/I(H)$ is isometrically isomorphic to $A_\Phi(H).$ Thus, by Corollary \ref{cor6}, $A_\Phi(H)$ is reflexive if and only if it is finite-dimensional, which, by Theorem \ref{thm8}, is equivalent to saying that $H$ is finite. Hence, the result follows.
	\end{proof}
	By choosing the open subgroup to be $G$ in the last corollary, we obtain the following.
	\begin{cor}\label{Ref_Chac}
		Let $G$ be an amenable group and let $\Phi$ satisfy the MA-condition. Then $A_\Phi(G)$ is reflexive if and only if $G$ is finite.
	\end{cor}
	
	\section{Arens regularity of $A_\Phi(G)$}
	
	In this section, we study the Arens regularity of $A_\Phi(G),$ in the spirit of \cite{For1}. We show that if $A_\Phi(G)$ is Arens regular, then $G$ is discrete. We also provide a necessary and sufficient condition for $A_\Phi(G)$ to be Arens regular in terms of the spaces $PM_\Psi(G)$ and $B_\Phi(G).$ 
	
	We begin this section by the following lemma which shows that many natural subspaces of $A_\Phi(G)$ are Arens regular, once we assume the Arens regularity of $A_\Phi(G).$
	
	\begin{lem}\label{lem1}
		Suppose that the Banach algebra $A_\Phi(G)$ is Arens regular. Then the following spaces are also Arens regular:
		\begin{enumerate}[(i)]
			\item Any closed ideal $I$ of $A_\Phi(G).$
			\item For any open subgroup $H$ of $G,$ the Banach algebra $A_\Phi(H).$
			\item For any compact normal subgroup $K$ of $G,$ the Banach algebra $A_\Phi(G/K).$
		\end{enumerate}
	\end{lem}
	\begin{proof}
		This follows directly from \cite[Pg. 312]{DH}, \cite[Lemma 4.1]{RLSK1} and \cite[Theorem 4.3]{RLSK1}.
	\end{proof}
	In the next result, we show that discreteness of the group $G$ is necessary for $A_\Phi(G)$ to be Arens regular.
	\begin{thm}\label{thm2}
		If $A_\Phi(G)$ is Arens regular, then $G$ is discrete.
	\end{thm}
	\begin{proof}
		First, let us assume that $G$ is second countable. Since $A_\Phi(G)$ is Arens regular, by \cite[Theorem 1]{DH}, for each $T \in PM_\Psi(G),$ the mapping $u \mapsto u\cdot T$ is weakly compact, which implies that $PM_\Psi(G) = WAP_\Psi(\widehat{G}).$ Therefore, $UCB_\Psi(\widehat{G}) \subseteq WAP_\Psi(\widehat{G})$ and hence, by \cite[Theorem 4.16]{RLSK3}, $G$ is discrete.
		
		Now, let $G$ be an arbitrary locally compact group. Let $G_0$ be a subgroup of $G$ which is $\sigma$-compact and open. By Lemma \ref{lem1}$\, (ii), \, A_\Phi(G_0)$ is Arens regular. If $G_0$ is second countable, then $G_0$ is discrete, and so is $G.$ If $G_0$ is not second countable, by \cite[Theorem 8.7]{HR1}, $G_0$ has a compact normal subgroup $N$ such that $m(N) = 0$ and $G_0/N$ is second countable. Since $A_\Phi(G_0)$ is Arens regular, again by Lemma \ref{lem1}$\, (i), \, A_\Phi(G_0/N)$ is also Arens regular. The preceding paragraph forces us to conclude that $G_0/N$ is discrete or, equivalently, $N$ is open, which contradicts that $m(N) = 0.$ Thus, $G$ is discrete.
	\end{proof}
	We now characterize discrete groups in terms of $A_\Phi(G).$
	\begin{thm}\label{thm3}
		The Banach algebra $A_\Phi(G)$ is an ideal in its second dual if and only if $G$ is discrete.
	\end{thm}
	\begin{proof}
		Suppose that $A_\Phi(G)$ is an ideal in $A_\Phi(G)^{**} = PM_\Psi(G)^*.$ Let $u \in A_\Phi(G)$ be such that $u(e) = 1 = \|u\|_{A_\Phi}.$ Consider the set $$K = \{m \odot u: m \, \, \text{is a mean on} \, \, PM_\Psi(G)\}.$$ It is then clear that $K \subseteq A_\Phi(G) \subseteq A_\Phi(G)^{**}.$ Since $\{m \in PM_\Psi(G)^*: m(I) = 1 = \|m\| \}$ is w$^*$-compact by Banach-Alaoglu theorem and the mapping $m \mapsto m \odot u$ is weak$^*$-weak$^*$-compact, the set $K$ is also w$^*$-compact. Note that the set $M_{A_\Phi}(G) = \{ v \in A_\Phi(G): v(e) = 1 = \|v\|_{A_\Phi}\}$ acts on $K$ as an abelian semigroup of weakly continuous affine operators $\{\varphi_v\}$ and since $K$ is convex, by Markov-Kakutani fixed point theorem, there exists $v_0 \in K$ such that $\varphi_v(v_0) = v_0$ for all $v \in M_{A_\Phi}(G).$ For $x \in G\setminus\{e\},$ choose $v_1 \in A_\Phi(G)$ such that $v_1 \in M_{A_\Phi}(G)$ and $v_1(x) = 0$ (see \cite[Prop. 5.5]{RLSK1}). Note that $v_0(x) = v_1(x) \, v_0(x) = 0.$ This helps us to conclude that $v_0 = \chi_{\{e\}}.$ Thus, $G$ is discrete.
		
		We shall now prove the converse. Suppose that $G$ is discrete. By \cite[Corollary 4.13]{RLSK3}, $UCB_\Psi(\widehat{G}) = PF_\Psi(G).$ Now, for $u \in A_\Phi(G), \, T \in PM_\Psi(G)$ and $m \in PM_\Psi(G)^*,$ $$\langle T,m \odot u \rangle = \langle u \cdot T,m \rangle.$$ Note that $u \cdot T \in UCB_\Psi(\widehat{G})$ and therefore, $m \in UCB_\Psi(\widehat{G})^* = PF_\Psi(G)^* = W_\Phi(G) \subseteq B_\Phi(G),$ where the last inclusion follows from \cite[Corollary 3.2]{RLSK4}. Thus, there exists $v \in B_\Phi(G)$ such that $\langle u \cdot T,m \rangle = \langle v,u \cdot T\rangle.$ Now $$\langle T,m \odot u \rangle = \langle u \cdot T,m \rangle = \langle v,u \cdot T \rangle = \langle uv,T \rangle = \langle T,uv \rangle,$$ which implies that $m \odot u = uv \in A_\Phi(G).$ Further, it is clear that $\langle T,m \odot u \rangle = \langle T,u \odot m \rangle,$ i.e., $A_\Phi(G)$ is in the centre of $PM_\Psi(G)^*.$ In particular, $A_\Phi(G)$ is a two-sided ideal in $PM_\Psi(G)^*.$
	\end{proof}
	\begin{thm}\label{thm4}
		Let $G$ be an amenable group and let $\Phi$ satisfy the MA-condition. Then $A_\Phi(G)$ is Arens regular if and only if $PM_\Psi(G)^*$ and $B_\Phi(G)$ are isomorphic as Banach algebras with equivalent norms.
	\end{thm}
	\begin{proof}
		Suppose that $A_\Phi(G)$ is Arens regular. By Theorem \ref{thm2}, $G$ is discrete and hence, by \cite[Corollary 4.13]{RLSK3}, $UCB_\Psi(\widehat{G}) = PF_\Psi(G).$ Since $G$ is amenable, by \cite[Corollary 4.13]{RLSK3}, $PF_\Psi(G) \subseteq WAP_\Psi(\widehat{G}) \subseteq UCB_\Psi(\widehat{G}).$ Further, using the Arens regularity, by \cite[Theorem 1]{DH}, $WAP_\Psi(\widehat{G}) = PM_\Psi(G).$ Putting all these together, we see that $$PF_\Psi(G) = WAP_\Psi(\widehat{G}) = UCB_\Psi(\widehat{G}) = PM_\Psi(G).$$ Now, by \cite[Corollary 3.3]{RLSK4}, $B_\Phi(G)$ and $W_\Phi(G) = PF_\Psi(G)^*$ are isomorphic as Banach spaces with equivalent norms. It is a routine to check, using the corresponding definitions, that $$\langle f,u \odot v \rangle = \langle f,u \cdot v \rangle \, \forall \, u,v \in B_\Phi(G) \, \text{and} \, f \in L^1(G).$$ As $L^1(G)$ is norm dense in $PF_\Psi(G)$ and $PF_\Psi(G) = PM_\Psi(G),$ we have $u \odot v = u \cdot v$ for all $u,v \in B_\Phi(G).$ Thus, $PM_\Psi(G)^*$ and $B_\Phi(G)$ are isomorphic as Banach algebras with equivalent norms. 
		
		The converse follows from \cite[Proposition 1]{DH}, hence proving the result.
	\end{proof}
	The next theorem provides a sufficient condition for the space $PM_\Psi(G)$ to have a unique topological invariant mean. As its proof follows verbatim to \cite[Theorem 3.12]{For1}, we omit the proof.
	\begin{thm}
		If $I$ is a non-zero ideal of $A_\Phi(G),$ which is also Arens regular, then $PM_\Psi(G)$ has a unique topological invariant mean.
	\end{thm}
	We now prove some corollaries by assuming the Arens regularity of $A_\Phi(G)$. The first one is about the Radon-Nikod\'ym property for $PM_\Psi(G).$
	\begin{cor}\label{cor4.6}
		Let $G$ be an amenable group for which $A_\Phi(G)$ is Arens regular. Then $PM_\Psi(G)$ has the Radon-Nikod\'ym property.
	\end{cor}
	\begin{proof}
		As $A_\Phi(G)$ is Arens regular, by Theorem \ref{thm2} and \ref{thm3}, it is a two-sided ideal in $PM_\Psi(G)^*.$ Also, as $G$ is amenable, by \cite[Theorem 3.1]{RLSK2}, $A_\Phi(G)$ has a bounded approximate identity. Now, the conclusion follows from \cite[Corollary 3.7]{ulg}.
	\end{proof}
	In the next corollary, we provide an equivalent condition for the space $A_\Phi(G)$ to be weakly sequentially complete.
	\begin{cor}
		Let $G$ be an amenable group for which $A_\Phi(G)$ is Arens regular and that $\Phi$ satisfies the MA-condition. Then $A_\Phi(G)$ is weakly sequentially complete if and only if $G$ is finite.
	\end{cor}
	\begin{proof}
		If $G$ is finite, then there is nothing to prove. Suppose that $A_\Phi(G)$ is weakly sequentially complete. Since $G$ is amenable, $A_\Phi(G)$ has an approximate identity $\{u_\alpha\}_{\alpha \in \Lambda},$ which is also bounded, by \cite[Theorem 3.1]{RLSK2}. For each $\alpha \in \Lambda,$ let $$A_\alpha = \{ u_\alpha \cdot T: T \in PM_\Psi(G), \|T\| \leq 1 \}.$$ It is clear that $A_\alpha$ is convex. Since $A_\Phi(G)$ is Arens regular, by Theorem \ref{thm3} and \cite[Lemma 3, Pg. 318]{DH}, $A_\alpha$ is weakly compact. Therefore, by \cite[Theorem 3.6.1, Pg. 60]{Bou}, each $A_\alpha$ has the RNP. Thus, by \cite[Corollary 3.9]{ulg}, the Banach algebra $A_\Phi(G)$ is reflexive. Hence, $G$ is finite by Corollary \ref{cor9}.
	\end{proof}
	Our final corollary is about the separability of $PM_\Psi(G).$
	\begin{cor}
		Let $G$ be a countable group, which is also amenable. If $A_\Phi(G)$ is Arens regular, then $PM_\Psi(G)$ is separable.
	\end{cor}
	\begin{proof}
		The assumption in the statement, along with Corollary \ref{cor4.6}, forces the conclusion that $PM_\Psi(G)$ has the Radon-Nikod\'ym property. The needed conclusion now follows from \cite[\textsection2]{ulg}.
	\end{proof}
	Before we move on to the last result of this section which characterizes the Arens regularity of $A_\Phi(G),$ here is a straightforward lemma.
	\begin{lem}\label{lem9}
		If $G$ is discrete, then $UCB_\Psi(\widehat{G}) = \overline{\text{span}\{\lambda_\Psi(x) : x \in G\}}.$
	\end{lem} 
	\begin{proof}
		The inclusion $UCB_\Psi(\widehat{G}) \subseteq \overline{\text{span}\{\lambda_\Psi(x) : x \in G\}}$ follows from Theorem \ref{thm3} and \cite[Theorem 4.1]{LU}.
		
		For the other way inclusion, it is enough to show that $\{\lambda_\Psi(x) : x \in G\} \subseteq UCB_\Psi(\widehat{G}).$ So, let $x \in G.$ Choose $u \in A_\Phi(G)$ such that $u(x) = 1.$ Then, for any $v \in A_\Phi(G),$ $$\langle v,\lambda_\Psi(x)\rangle = v(x) = v(x) \, u(x) = \langle uv,\lambda_\Psi(x) \rangle = \langle v,u \cdot \lambda_\Psi(x)\rangle.$$ Then, $\lambda_\Psi(x) = u \cdot \lambda_\Psi(x) \in UCB_\Psi(\widehat{G}).$ Hence the proof.
	\end{proof}
	Here is the characterization result. It includes Theorem \ref{thm4}. This result is the Orlicz analogue of \cite[Corollary 8.3]{LU}.
	\begin{thm}
		Let $G$ be an amenable group. Then TFAE:
		\begin{enumerate}[(i)]
			\item The Banach algebra $A_\Phi(G)$ is Arens regular.
			\item The spaces $PM_\Psi(G)$ and $PF_\Psi(G)$ are equal.
			\item The spaces $PM_\Psi(G)$ and $\overline{\text{span}\{\lambda_\Psi(x) : x \in G\}}$ are equal.
		\end{enumerate}
	\end{thm}
	\begin{proof}
		$(i) \iff (ii).$ Let $A_\Phi(G)$ be Arens regular. Then, the proof of $(ii)$ is already contained in Theorem \ref{thm4}.
		
		Conversely, let condition $(ii)$ hold. Then, by \cite[Corollary 4.13]{RLSK3}, we have $PM_\Psi(G) = WAP_\Psi(\widehat{G}).$ Therefore, by \cite[Theorem 1]{DH}, $A_\Phi(G)$ is Arens regular.
		
		$(i) \iff (iii).$ Let $A_\Phi(G)$ be Arens regular. Then, by Theorem \ref{thm2}, $G$ is discrete. The conclusion now follows from Lemma \ref{lem9} and the equivalence between $(i)$ and $(ii).$
		
		Conversely, let condition $(iii)$ hold. Then, in order to show that $A_\Phi(G)$ is Arens regular, by \cite[Theorem 1]{DH}, it is enough to show that for each $x \in G, \, \lambda_\Psi(x) \in AP_\Psi(\widehat{G}),$ i.e., the map $u \mapsto u \cdot \lambda_\Psi(x)$ from $A_\Phi(G)$ to $PM_\Psi(G)$ is a compact operator. Now, for any $u,v \in A_\Phi(G)$ and $x \in G,$ $$\langle v,u \cdot \lambda_\Psi(x) \rangle = \langle uv,\lambda_\Psi(x)\rangle = u(x) \, v(x) = u(x) \langle v,\lambda_\Psi(x)\rangle = \langle v, \langle u,\lambda_\Psi(x)\rangle \lambda_\Psi(x)\rangle.$$ Therefore, for any $u \in A_\Phi(G), \, u \cdot \lambda_\Psi(x) = \langle u,\lambda_\Psi(x)\rangle \lambda_\Psi(x).$ In particular, the mapping $u \mapsto u \cdot \lambda_\Psi(x)$ is a bounded operator with rank one, thereby proving its compactness.
	\end{proof}
	
	\section{Arens regularity on discrete groups}
	
	Our main aim in this section is to give some sufficient conditions
	so that $A_\Phi(G)$ is Arens regular. We shall begin with a simple result. This can be considered as a converse to Lemma \ref{lem1}.
	\begin{lem}
		The Banach algebra $A_\Phi(G)$ is Arens regular if and only if $G$ is discrete and for each countable subgroup $H$ of $G,$ $A_\Phi(H)$ is Arens regular.
	\end{lem}
	\begin{proof}
		We shall prove only the backward implication, as the forward is already proved in Lemma \ref{lem1} and Theorem \ref{thm2}. Now, to prove the backward, our strategy is to make use of \cite[Theorem 1]{DH}.
		
		Let $\{u_n\}$ and $\{v_n\}$ be two sequences in $A_\Phi(G) \cap \mathbf{C}_c(G)$ and let $T \in PM_\Psi(G).$ Now, let $K = \bigg( \bigcup\limits_n \text{supp}(u_n) \bigg) \bigcup \bigg( \bigcup\limits_m \text{supp}(v_m) \bigg).$ Since $G$ is discrete, $K$ is countable. Let $H$ denote the subgroup generated by $K.$ Then $H$ is open, and hence, by assumption, $A_\Phi(H)$ is Arens regular. Therefore, by \cite[Theorem 1]{DH}, $$\lim\limits_{n} \lim\limits_{m} T(u_n \, v_m) = \lim\limits_{m} \lim\limits_{n} T(u_n \, v_m).$$ Hence, the proof follows.
	\end{proof}
	From now onwards, we assume throughout this section that the group $G$ is discrete.
	\begin{rem}
		Since $G$ is discrete, $l^1(G)$ naturally sits inside $A_\Phi(G).$ In fact, if $f \in l^1(G),$ then $\sum\limits_{x \in G} | f(x) | < \infty,$ i.e., there exists a countable set $E \subseteq G$ such that $f(x) = 0$ for all $x \notin E.$ Thus, if we assume that $E = \{ x_n: n \in \mathbb{N}\},$ then $f(x) = \sum\limits_{n \in \mathbb{N}} f(x_n) \delta_{x_n}(x),$ where $$\delta_{x_n}(x) = 
		\begin{cases} 
			1 & \text{if} \, \, x = x_n \\
			0 & \text{otherwise.} 
		\end{cases}$$ 
		Observe that $\delta_{x_n}(x) = \chi_{\{x_n\}} \ast \check{\chi}_{\{e\}}(x).$ Therefore, $$f = \sum\limits_{n \in \mathbb{N}} f(x_n) \chi_{\{x_n\}} \ast \check{\chi}_{\{e\}}.$$ Also,
		\begin{align*}
			\|f\|_{A_\Phi} &\leq \sum\limits_{n \in \mathbb{N}} |f(x_n)| \, N_\Phi(\chi_{\{x_n\}}) \, \|\chi_{\{e\}}\|_{\Psi} \\
			&= \sum\limits_{n \in \mathbb{N}} |f(x_n)| \, [\Phi^{-1}(1)]^{-1} \, [\Phi^{-1}(1)] = \sum\limits_{n \in \mathbb{N}} |f(x_n)| = \|f\|_1 < \infty.
		\end{align*}
		This implies that $f \in A_\Phi(G)$ and $\|f\|_{A_\phi} \leq \|f\|_1,$ i.e., $l^1(G) \subseteq A_\Phi(G).$
	\end{rem}
	\begin{rem}\label{rem5.3}
		We remark here that $(l^1(G),\cdot,\|\cdot\|_{A_\Phi})$ is only a normed algebra but not a Banach algebra. Here, $`\cdot$' denotes the pointwise multiplication. Since $l^1(G)$ is norm dense in $A_\Phi(G),$ it follows that $(l^1(G),\cdot,\|\cdot\|_{A_\Phi})$ is Arens regular if and only if $A_\Phi(G)$ is Arens regular.
	\end{rem}
	Motivated by the above remark, we first consider the Arens regularity of $(l^1(G),\cdot,\|\cdot\|_1).$
	\begin{lem}\label{lem5.4}
		The Banach algebra $(l^1(G),\cdot,\|\cdot\|_1)$ is Arens regular.
	\end{lem}
	\begin{proof}
		Let $\{f_n\}$ and $\{g_n\}$ be any two sequences from the unit ball of $l^1(G)$ and let $h \in l^\infty(G)$ be such that both the limits $$\lim\limits_n \lim\limits_m \sum\limits_{x \in G} h(x) f_n(x) g_m(x) \, \text{and} \, \lim\limits_m \lim\limits_n \sum\limits_{x \in G} h(x) f_n(x) g_m(x)$$ exists. As $l^1(G) \cong C_0(G)^*,$ the unit ball of $l^1(G)$ is w$^*$-compact and hence, the sequences $\{f_n\}$ and $\{g_n\}$ have subnets $\{f_{n_\alpha}\}$ and $\{g_{n_\beta}\}$ such that they converge in w$^*$-topology to, say $f$ and $g,$ respectively. Now, $$\lim\limits_n \lim\limits_m \sum\limits_{x \in G} h(x) f_n(x) g_m(x) = \lim\limits_n \sum\limits_{x \in G} h(x) f_n(x) g(x) = \sum\limits_{x \in G} h(x) f(x) g(x).$$ Similarly, $\lim\limits_m \lim\limits_n \sum\limits_{x \in G} h(x) f_n(x) g_m(x) = \sum\limits_{x \in G} h(x) f(x) g(x).$ Now, the conclusion follows from \cite[Theorem 1]{DH}.
	\end{proof}
	\begin{nota}
		Let $B_{\Phi,1}$ denote the unit ball of $(l^1(G),\cdot,\|\cdot\|_{A_\Phi}).$
	\end{nota}
	We observed that $l^1(G)\subseteq A_\Phi(G).$ Consequently, we first provide a condition under which $(l^1(G),\cdot,\|\cdot\|_{A_\Phi})$ is Arens regular.
	\begin{thm}\label{thm5.5}
		Suppose that $B_{\Phi,1}$ is bounded w.r.t. $\|\cdot\|_1,$ then $(l^1(G),\cdot,\|\cdot\|_{A_\Phi})$ is Arens regular.
	\end{thm}
	\begin{proof}
		Let $f \in l^1(G)$ and $T \in (l^1(G),\cdot,\|\cdot\|_{A_\Phi})^*.$ Then, $$|\langle f,T \rangle | \leq \|T\| \, \|f\|_{A_\Phi} \leq \|T\| \|f\|_1,$$ i.e., $T \in (l^1)^* \cong l^\infty.$ Now, let $\{f_n\}$ and $\{g_n\}$ be any two sequences from $B_{\Phi,1}.$ As $B_{\Phi,1}$ is bounded w.r.t. $\|\cdot\|_1,$ the sequences $\{f_n\}$ and $\{g_n\}$ are bounded in $(l^1(G),\cdot,\|\cdot\|_1).$ By Lemma \ref{lem5.4}, $(l^1(G),\cdot,\|\cdot\|_1)$ is Arens regular and hence, by \cite[Theorem 1]{DH}, $$\lim\limits_n \lim\limits_m T(f_n \, g_m) = \lim\limits_m \lim\limits_n T(f_n \, g_m).$$ This implies that $(l^1(G),\cdot,\|\cdot\|_{A_\Phi})$ is Arens regular, again by \cite[Theorem 1]{DH}.
	\end{proof}
	By combining Remark \ref{rem5.3} and Theorem \ref{thm5.5} we obtain the following corollary.
	\begin{cor}
		Suppose that $B_{\Phi,1}$ is bounded w.r.t. $\|\cdot\|_1.$ Then $A_\Phi(G)$ is Arens regular.
	\end{cor}
	
	\section{Arens regularity of $PF_\Psi(G)$ and related results}
	
	In the final section of this paper, we consider the algebra $PF_\Psi(G).$ We shall prove several inclusion results for this algebra. On the other hand, we shall also see a similar result related to the algebra $PM_\Psi(G).$ The results of this section are Orlicz analogues of the results given in \cite{LU}.
	
	We shall begin with a lemma.
	
	\begin{lem}\label{lem5.1}
		For any countable group $G,$ the containment $$wap(\mathbf{B}(L^\Psi(G))) \supseteq L^\Phi(G) \, \widehat{\otimes} \, L^\Psi(G)$$ holds.
	\end{lem}
	\begin{proof}
		Since $L^\Psi(G)$ is reflexive, by \cite[Theorem 2.6.23]{Da}, the algebra $\mathbf{K}(L^\Psi(G))$ is Arens regular. Further, by assumption, $G$ is countable; hence, by \cite{LT}, $L^\Psi(G)$ has the approximation property. Thus, the proof follows from \cite[Pg. 195]{DF} and the fact that $\mathbf{K}(L^\Psi(G))^* \subseteq wap(\mathbf{B}(L^\Psi(G))).$
	\end{proof}
	As a consequence, we have the following corollary.
	\begin{cor}\label{APhisubsetwap}
		For a countable group $G, A_\Phi(G) \subseteq wap(PF_\Psi(G)).$
	\end{cor}
	\begin{proof}
		This follows from Lemma \ref{lem5.1} and the isometric isomorphism between the Banach space $A_\Phi(G)$ and the quotient space $(L^\Phi(G) \, \widehat{\otimes} \, L^\Psi(G))/\text{ker}(\varphi)$ where $\varphi : L^\Phi(G) \, \widehat{\otimes} \, L^\Psi(G) \to A_\Phi(G),$ given by $\varphi(\sum f_n \otimes g_n) = \sum f_n \ast \check{g_n}.$
	\end{proof}
	It is a known fact that for a compact group $G$ and $f \in L^1(G),$ the convolution operator $\lambda_\Psi(f)$ is compact. Consequently, we obtain the following theorem which states that the multiplication map on $PF_\Psi(G)$ is compact.
	\begin{thm}\label{thm5.3}
		Let $G$ be a compact group. Then, for each $T \in PF_\Psi(G),$ the mapping $S \mapsto S \circ T$ (also, $S \mapsto T \circ S$) from $PF_\Psi(G)$ to $PF_\Psi(G)$ is compact.
	\end{thm}
	\begin{proof}
		As the proof of this theorem follows along the same lines of \cite[Lemma 8.6]{LU}, we omit the proof. 
	\end{proof}
	Consequently, we have the following inclusion result.
	\begin{cor}\label{cor5.4}
		Let $G$ be a compact group. Then $ap(PF_\Psi(G)) \subseteq W_\Phi(G).$
	\end{cor}
	\begin{proof}
		Since $L^1(G)$ has a bounded approximate identity and as $L^1(G)$ is dense in $PF_\Psi(G),$ the algebra $PF_\Psi(G)$ also has a bounded approximate identity. By Theorem \ref{thm5.3}, for each $T \in PF_\Psi(G),$ the sets $\{S \circ T : S \in PF_\Psi(G) \, , \, \|S\| \leq 1\}$ and $\{T \circ S : S \in PF_\Psi(G) \, , \, \|S\| \leq 1\}$ are weakly compact. Thus, $PF_\Psi(G)$ is a two-sided ideal in its second dual, and therefore, by \cite[Corollary 3.2]{ulg}, the conclusion follows.
	\end{proof}
	Our next result is about the Arens regularity of $PF_\Psi(G).$
	\begin{cor}
		Let $G$ be a compact group. Then $PF_\Psi(G)$ is Arens regular and $W_\Phi(G)$ has the Radon-Nikod\'ym property. 
	\end{cor}
	\begin{proof}
		Since $G$ is compact, for each $f \in L^1(G),$ the mapping $\lambda_\Psi(f),$ given by $\lambda_\Psi(f)(g) = f \ast g,$ is a compact operator. Hence, $PF_\Psi(G) \subseteq \mathbf{K}(L^\Psi(G)).$ Thus, the Arens regularity of $PF_\Psi(G)$ is a consequence of \cite[Theorem 2.6.23]{Da} and the fact that $L^\Psi(G)$ is reflexive. Also, from the proof of Corollary \ref{cor5.4}, it follows that $PF_\Psi(G)$ is a two-sided ideal in $W_\Phi(G)^*$ and thus, by \cite[Corollary 3.7]{ulg}, $W_\Phi(G)$ has the RNP.
	\end{proof}
	Here is the final result of this paper. This is the $PM_\Psi(G)$ analogue of Corollary \ref{APhisubsetwap}. However, the inclusion result here actually characterizes compact groups. This is the Orlicz analogue of \cite[Theorem 8.8]{LU}.
	\begin{thm}
		Let $G$ be a locally compact group. Then TFAE:
		\begin{enumerate}[(i)]
			\item $G$ is compact.
			\item $A_\Phi(G) \subseteq ap(PM_\Psi(G)).$
			\item The multiplication is jointly continuous in the unit ball of $PM_\Psi(G)$ equipped with the w$^*$-topology.
		\end{enumerate}
	\end{thm}
	\begin{proof}
		$(i) \implies (ii).$ In order to prove $(ii),$ it is enough to prove that $u \in ap(PM_\Psi(G)),$ when the function $u \in A_\Phi(G)$ is of the form $u = f \ast \check{g},$ where $f \in L^\Phi(G)$ and $g \in L^\Psi(G).$
		
		For $S,T \in PM_\Psi(G),$
		\begin{align*}
			\langle S, J(u) \cdot T \rangle &= \langle S \circ T,J(u) \rangle = \langle u,S \circ T \rangle = \langle f \ast \check{g},S \circ T \rangle \\
			& = \langle S \circ T(g),f \rangle = \langle f \ast \check{T(g)},S \rangle = \langle v,S \rangle = \langle S,J(v) \rangle,
		\end{align*}
		where $v = f \ast \check{T(g)} \in A_\Phi(G)$ and $J$ is the canonical map from $A_\Phi(G)$ to $PM_\Psi(G)^*.$ Thus, $J(u) \cdot T = J(v) \in PM_\Psi(G)^*$ for some $v \in A_\Phi(G).$ In fact, $J(u) \cdot T \in J(A_\Phi(G)).$
		
		Let $h \in L^\Psi(G).$ Since $G$ is compact, $h \in L^1(G)$ and hence, $\lambda_\Psi(h) \in PM_\Psi(G).$ Let $w = J(u) \cdot \lambda_\Psi(h).$ It can be shown as above that, for $T \in PM_\Psi(G),$ $w \cdot T \in J(A_\Phi(G)).$ We now claim that $w \in ap(PM_\Psi(G)),$ i.e., the set $\{w \cdot T: T \in PM_\Psi(G), \|T\| \leq 1\}$ is relatively compact in $J(A_\Phi(G)).$ So, let $S \in PM_\Psi(G).$ Then,
		\begin{align*}
			S(w \cdot T) &= \langle S,w \cdot T \rangle = \langle S\circ T,w \rangle = \langle S\circ T,J(u) \cdot \lambda_\Psi(h) \rangle \\
			&= \langle u, S \circ T \circ \lambda_\Psi(h) \rangle =
			\langle f \ast \check{g},S \circ T \circ \lambda_\Psi(h) \rangle \\
			&= \langle S \circ T \circ \lambda_\Psi(h)(g),f \rangle = \langle S \circ T(h \ast g),f \rangle \\
			&= \langle f \otimes (T(h \ast g)),S \rangle = \langle f \otimes (T(h) \ast g),S \rangle = S(f \otimes (T(h) \ast g)).
		\end{align*}
		Thus, it is enough to prove that if $\{S_\alpha\}$ is a net from the unit ball of $PM_\Psi(G)$ that converges to $0$ in the w$^*$-topology, then $\{S_\alpha\}$ converges uniformly to $0$ on $\{w \cdot T: T \in PM_\Psi(G), \|T\| \leq 1\}.$ Note that, for $h,g \in L^\Psi(G),$ the mapping $T \mapsto T(h) \ast g$ is compact and hence, the set $\{T(h) \ast g: T \in PM_\Psi(G), \|T\| \leq 1\}$ is relatively compact in $L^\Psi(G).$ This implies that the set $\{f \otimes (T(h) \ast g): T \in PM_\Psi(G), \|T\| \leq 1\}$ is relatively compact in $L^\Phi(G) \otimes L^\Psi(G).$ Since $S_\alpha \xrightarrow{w^*} 0$ in $\mathbf{B}(L^\Psi(G)),$ we have, $\lim\limits_\alpha \sup\limits_{\|T\| \leq 1} S_\alpha(f \otimes (T(h) \ast g)) = 0.$ From the above computations, it follows that $\lim\limits_\alpha \sup\limits_{\|T\| \leq 1} S_\alpha(w \cdot T) = 0$ and hence, our claim follows.
		
		Now, let $\{h_\alpha\}$ be an approximate identity in $L^\Psi(G).$ Then, for $T \in PM_\Psi(G),$
		\begin{align*}
			\langle T,u \cdot \lambda_\Psi(h_\alpha) \rangle &= \langle u,T \circ \lambda_\Psi(h_\alpha) \rangle = \langle f \ast \check{g},T \circ \lambda_\Psi(h_\alpha) \rangle \\
			&= \langle T \circ \lambda_\Psi(h_\alpha)(g), f\rangle = \langle T(h_\alpha \ast g),f \rangle \rightarrow \langle T(g),f \rangle = \langle T,u \rangle,
		\end{align*}
		i.e., $u \cdot \lambda_\Psi(h_\alpha) \to u$ weakly in $A_\Phi(G).$ Since, $ap(PM_\Psi(G))$ is a closed subspace of $PM_\Psi(G)^*,$ the conclusion follows.
		
		The rest of the proof follows along the same lines of \cite[Theorem 8.8]{LU}. Hence, we omit it.
	\end{proof}
	
	\section*{Data Availability} 
	Data sharing does not apply to this article as no datasets were generated or analyzed during the current study.
	
	\section*{competing interests}
	The authors declare that they have no competing interests.

\end{document}